\documentclass[10pt,a4paper,reqno]{amsart}
\usepackage{amsfonts,amsthm,latexsym,amsmath,amssymb,amscd,amsmath, epsf}
\usepackage{graphicx, epigraph}

\usepackage{graphicx}
\newtheorem{theorem}{Theorem}
\newtheorem{proposition}[theorem]{Proposition}
\newtheorem{lemma}[theorem]{Lemma}
\newtheorem{definition}{Definition}
\newtheorem{corollary}{Corollary}

\newtheorem{conjecture}{Conjecture}

\newtheorem{problem}{Problem}

\newcommand{\bR}{\mathbb{R}}
\newcommand{\bC}{\mathbb{C}}
\newcommand{\be}{\begin{equation}}
\newcommand{\ee}{\end{equation}}
\newcommand {\al}{\alpha}

\newcommand{\Log}{\operatorname{Log}}

\renewcommand{\mod}{\mathop{\rm \ mod}}

\newcommand {\si} {\sigma}
\newcommand{\bsi}{{\bar \si}}
\newcommand{\bZ}{\mathbb Z}

\begin{document}
          \numberwithin{equation}{section}

          \title[ Problems around  polynomials - the good, the bad and the ugly\dots  ]
          {Problems around  polynomials - the good, the bad and the ugly\dots   }

\author[B. Shapiro]{Boris Shapiro}
\address{   Department of Mathematics,
            Stockholm University,
            S-10691, Stockholm, Sweden}
\email{shapiro@math.su.se}

%\begin{abstract}   In what follows I formulate a number of problems related to (mainly univariate) polynomials which  are   %motivated by my research of the last decade. 

%\end{abstract}
\maketitle

\dedicatory{The Russian style of formulating mathematical problems  means that nobody will be able to simplify your formulation as opposed to  the French style which means that nobody will be able to generalize it, -- Vladimir Arnold\footnote {The title of the paper alluding to the highly recommended spaghetti western by S.~Leone, reflects the very personal taste of the author concerning mathematical problems which should be taken with a grain of salt. In his opinion  a "good`` problem  has an elegant short formulation and is (hopefully) solvable, while a "bad`` one looks equally stimulating but seems unaccessible at the moment. Finally an "ugly`` problem apparently has a non inspiringly complicated answer which can hardly stumulate a further development. But one can not really tell  untill the problem itself has been actually solved! But in any case,  the esthetic feeling about mathematical problems and their solutions has to be taken seriously. }  }

\section{Problems}  

\noindent
{\bf I. Around Maxwell's conjecture.}  In Section 133 of \cite{Maxwell},  J.~C.~Maxwell formulated the following claim, but provided it will an incomplete proof (see details in \cite{GNSh}).

\begin{conjecture}[ Maxwell, seems bad, no tools]\label{conj:Maxw}
For any system of $N$ isolated fixed point charges in $\bR^3$, the number of  points of equilibrium (assumed finite) of the created electrostatic field
$$E(\bar x)=\sum_{i=1}^N\frac{\xi_i(\bar x-\bar x_i)}{|\bar x-\bar x_i|^3},$$
 is at most $(N-1)^2$. Here $\xi_i$ is a charge placed at $\bar x_i\in \bR^3$ and $\bar x\in \bR^3$ is a variable vector.  
\end{conjecture}

Some very crude estimate on the maximal number of points of equilibrium is obtained in \cite{GNSh} using Khovanskii's fewnomial theory, see \cite{Kh}. 
Conjecture~\ref{conj:Maxw} is not even settled  for 3 positive point charges in which case all points of equilibrium  lie in the plane spanned by these charges. 
(Some special cases of 3 charges are settled in the recent literature, see \cite{Ki},\cite{Pe},\cite{Wa}.) 

 Observe that for charges of different signs in $\bR^3,$ the set of  points of equilibrium might be a space curve. The simplest example of that kind is a 4-tuple of charges placed at the vertices of a square in the $z$-plane with coordinates $(\pm 1, \pm 1, 0)$. If  one places the unit positive charges at $(1,1,0)$ and $(-1,-1,0)$ and the unit negative charges at 2 remaining corners, then the set of  points of equilibrium coincides with the $z$-axis. The following naive-looking question is not settled either.

\begin{conjecture}[folklore, very irritating]\label{conj:positive}
For any set of charges of the same sign in $\bR^n$, the set of its  points of equilibrium is finite. 
\end{conjecture}

The next claim is a very special one-dimensional case of a rather general Conjecture 1.7. of \cite{GNSh} which generalizes the above Conjecture~\ref{conj:Maxw} in many ways.

\begin{conjecture}[A.~Gabrielov, D.~Novikov,  B.~Sh., seems good, but no progress]\label{1-dim}
Let $(x_1,y_1),(x_2,y_2),\dots, (x_N,y_N)$ be a collection of points in $\bR^2$, $\xi_1, \xi_2, \dots , \xi_N$ be arbitrary real charges and $\al\ge 1/2$. Then  the rational  univariate function 
$$\Psi(x)=\sum_{i=1}^N\frac{\xi_i}{((x-x_i)^2+y_i^2)^{\al}},\; x\in \bR,$$
has at most $N$ local maxima on the whole real line. 

\end{conjecture} 

Observe that also this conjecture is not settled  in the simplest case $N=3$, $\al=1$ and all unit charges. The author has   an overwhelming  numerical evidence supporting the latter conjecture but no proof. 

\bigskip
\noindent
{\bf II. On real zeros of exponential sums.}  Consider the space of linear ordinary homogeneous differential equations with constant coefficients of order $k$, i.e.
\begin{equation}\label{eq:triv}
y^{(k)}+a_1y^{(k-1)}+\dots + a_ky=0
\end{equation}
where $a_1,\dots, a_k$ are arbitrary complex numbers. The next statement easily follows from the standard facts about the asymptotic zero distribution of exponential sums, see e.g., \cite{La}. 

\begin{lemma}\label{lm:simple} Every non-trivial solution of \eqref{eq:triv} has  finitely many (probably none) real zeros if and only if 
any two distinct characteristic roots of \eqref{eq:triv} have distinct real parts. 
\end{lemma}

Denote by $\Omega_k$ the set of all \eqref{eq:triv} satisfying the conditions of Lemma~\ref{lm:simple}.  (It is an open dense subset of $\bC^k$ with coordinates $(a_1,\dots , a_k)$.) 

\begin{problem}[B.~Sh., looks bad, but very important] 
Does there exist an upper bound for the number of real roots valid for all non-trivial solutions of all equations~\eqref{eq:triv}  in $\Omega_k$?
\end{problem}

The latter problem is open already for $k=3$. Observe that there exists a highly non-trivial and apparently far from being sharp  upper bound for the number of integer zeros of exponential polynomials obtained in \cite{Sch}.

\bigskip
\noindent
{\bf III. On isolated zeros of non-negative polynomials and sums of squares.} 

\begin{problem}[D.~Khavinson, I.~Itenberg,  B.~Sh., apparently bad]
 Find  the maximal possible number $\sharp(2k, l)$ of isolated zeros for real non-negative polynomials of degree $2k$ in $l$  variables.
\end{problem}

\begin{problem}[G.~Ottaviani,  B.~Sh., seems good]
 Find  the maximal possible number $\widetilde\sharp(2k, l)$ of isolated zeros for real non-negative polynomials of degree $2k$ in $l$  variables which are representable as the sums of squares of real polynomials of degree at most $k$. 
\end{problem}

A trivial observation is that $k^l\le \widetilde\sharp(2k,l)\le  \sharp(2k, l)\le (2k -1)^l.$  In the special case $l=2$ both problems were considered in an intruguing paper \cite{CLR} where  it was proven that $\widetilde \sharp(2k,2)=k^2$ and $\sharp(2k,2)\le \frac{3k(k-1)}{2}+1$. The latter inequality is obtained with the help of Petrovskii-Oleinik's inequality, see \cite{OlPe}. 

\begin{conjecture}[G.~Ottaviani,  B.~Sh., seems good]
For any number of variables, $\widetilde\sharp(2k,l)=k^l$. 
\end{conjecture} 

On the other hand, it seems difficult to determine the coefficient of  $k^l$ of the leading asymptotic term for $\sharp(2k,l),$ when $l$ is fixed and $k\to\infty$. For example, I doubt that $\sharp(2k,2)$ grows asymptotically as ${3k^2}/{2},$ when $k\to\infty$.

\bigskip
\noindent
{\bf IV. Hermite-Biehler problem.}
The well-known Hermite-Biehler theorem claims that a univariate monic polynomial $s$ of degree $k$ has all roots in the open upper half-plane if and only if  $s=p+iq.$ Here  $p$ and $q$ are real polynomials of degree $k$ and $k-1$ respectively with all real, simple and interlacing roots,  and $q$ has a negative leading coefficient.

\medskip
\noindent
\begin{problem}[S.~Fisk, seems bad, see \cite{Fi}, p.~575] Given a pair of real polynomials $(p,q),$ give restrictions on the location of the roots of $p+iq$ in terms of the location of the roots of $p$ and $q$.
\end{problem}

An example of such results can be found in \cite{KST}. Closely related important questions are:  (i) restrictions on the location of (complex) roots of the Wronskian $W(p,q);$  (ii)  description of the real univalent disks for the real rational function 
$\frac{p}{q}$.

\bigskip
\noindent
{\bf V. Mesh-related questions.} 
 By the {\bf mesh} of a  polynomial $p(x)$ with all real simple zeros we mean the minimal distance between its consecutive roots. 

\begin{conjecture}[P.~Br\"anden, I.~Krasikov, B.~Sh., hopefully good, see \cite{BKSh}]\label{pr:finitedegreeDif} A difference operator $T(p(x))=a_0p(x)+a_1p(x-1)+\cdots+a_kp(x-k)$ with constant coefficients preserves the set of  real-rooted polynomials of degree at most $m$ whose mesh is at least $1$  if and only if the polynomial $T((x)_m)$ is real-rooted and has mesh at least one.  Here $(x)_m=x(x-1)(x-2)\dots (x-m+1)$ is the $m$-th Pochhammer polynomial.  
\end{conjecture}
There is an alternative formulation of Conjecture \ref{pr:finitedegreeDif} which is maybe more attractive. Let 
$\nabla p(x)= p(x+1)-p(x)$ be the forward difference operator, and consider the following product on the space of polynomials of degree at most $d$:
$$
(p \bullet q)(x) = \sum_{k=0}^d (\nabla^kp)(0) \cdot(\nabla^{d-k}q)(x).
$$
Conjecture \ref{pr:finitedegreeDif} is equivalent to
\begin{conjecture}\label{pr:finitedegreeDif2} If $p$ and $q$ are real-rooted polynomials of degree at most $d$ and of mesh $\geq 1$, then so is $p \bullet q$. 
\end{conjecture}

\bigskip
\noindent
{\bf VI. Topology of the space of polynomials.}
 The famous Descartes' rule of signs claims that the number of positive roots of a real univariate polynomial does not exceed the number of sign changes in its sequence of coefficients. For simplicity let us only consider polynomials with all non-vanishing coefficients. An arbitrary ordered sequence 
$\bsi=(\si_0,\si_1,...,\si_d)$ of $\pm$-signs  is called a {\em sign pattern}. Given a sign pattern ${\bsi}$ as above, we call by its {\em Descartes' pair} $(p_\bsi,n_\bsi)$ the pair of non-negative integers counting sign changes and sign preservations of $\bsi$.  The Descartes' pair of $\bsi$ gives the upper bounds on the number of positive and negative roots   of any polynomial of degree $d$ whose signs of coefficients are given by   $\bsi$. (Observe that, for any $\bsi$, 
$p_\bsi+ n_\bsi=d$.) To any polynomial $q(x)$ with the sign pattern $\bsi,$ we associate  the pair $(pos_q,neg_q)$ giving the numbers of its positive and negative roots (counted multiplicities). Obviously $(pos_q,neg_q)$  satisfies the standard restrictions
 \begin{equation}\label{stand} 
pos_q\le p_\bsi,\;  pos_q\equiv p_\bsi (mod\, 2),\; neg_q\le n_\bsi,\; neg_q\equiv n_\bsi (mod\, 2).
\end{equation}

We call pairs $(pos, neg)$ satisfying \eqref{stand} {\em admissible} for $\bsi$.  
 It turns out that not for every pattern $\bsi,$ all its admissible  pairs  $(pos,neg)$  are realizable by polynomials with the sign pattern $\bsi$. 

\begin{problem}[seems bad, but might be ugly]
For a given sign pattern $\bsi,$  which admissible pairs $(pos,neg)$ are realizable by polynomials whose signs of coefficients are given by  $\bsi$?  
\end{problem}

The first non-realizable combination of a sign pattern and a pair $(pos,neg)$ occurs in degree $4$, see \cite{Gr}. Namely,  (up to the standard $\bZ_2\times \bZ_2$-action on the set of all sign patterns) the only non-realizable combination is $\bsi=(+,+,-,+,+)$ with the pair $(2,0)$. Based on our computer-aided results up to degree $10$, we can formulate the following claim. 

\begin{conjecture}[J.~Forsg\aa rd, V.~Kostov, B.~Sh,  hopefully good, see \cite{FKS}] \label{conj:main} For an arbitrary sign pattern $\bsi$, the only type of pairs $(pos,neg)$ which can be non-realizable has either $pos$ or $neg$ vanishing. In other words, for any sign pattern $\bsi$, each pair $(pos,neg)$ satisfying \eqref{stand}  with positive $pos$ and $neg$ is realizable.
\end{conjecture}

\bigskip
\noindent
{\bf VII. Tropical geometry.} % See \cite{FSh}. 
%\subsection{Tropical bound on the number of real roots} 

\begin{conjecture}[J.~Forsg\aa rd, B.~Sh., seems good,  see \cite{FSh}] 
 Let $f(z) = \sum_{k=0}^n a_k z^k$ be a polynomial with positive coefficients, and consider the related (weighted) 
 tropical polynomial 
 \[
  f_{trop}(x) =\max_{k}\left(\Log(a_k)+k x + \Log{n\choose k}\right).
 \]
 Then the number of real zeros of $f(z)$ does not exceed the number of points in the tropical variety defined by  $f_{trop}$, i.e. the number of corners of the continuous piecewise-linear function 
 $ f_{trop}(x),\; x\in \bR$.
\end{conjecture}

%\subsection{Descartes--Newton bound on the number of real roots}
\begin{conjecture}[J.~Forsg\aa rd, B.~Sh., seems good,  see \cite{FSh}] 
\label{ConjDiffInduktion}
 Let $f(z) = \sum_{k=0}^n a_k z^k$ be a polynomial with positive coefficients. Consider the differences
 \[
  \tilde c_k = (k+1)a_k^2 - k a_{k-1}a_{k+1},
 \]
 where $a_{-1} = a_{n+1} = 0$. Let $0=k_1 < k_2 < \dots < k_m = n$ be the sequence of indices such that $\tilde c_{k_i}$ 
 is positive, and let $v(f)$ be the number of changes in the sequence $\{k_i \mod 2\}_{i=0}^m$. Then the
 number of real zeros of $f(z)$ does not exceed $v(f)$.
\end{conjecture}

\begin{conjecture}[J.~Forsg\aa rd, B.~Sh., seems good,  see \cite{FSh}]] 
 Let $f(z) = \sum_{k=0}^n a_k z^k$ be a polynomial with positive coefficients. Consider the differences
 \[
  c_k = a_k^2 - a_{k-1}a_{k+1},
 \]
 where $a_{-1} = a_{n+1} = 0$. Let $0=k_1 < k_2 < \dots < k_m = n$ be the sequence of indices such that $c_{k_i}$ 
 is non-negative, and let $v(f)$ be the number of changes in the sequence $\{k_i \mod 2\}_{i=0}^m$. Then the
 number of real  zeros of $f(z)$ does not exceed $v(f)$.
\end{conjecture}

Observe that for polynomials with all positive coeffficients,  their real roots are negative.

\bigskip
\noindent
{\bf VIII. Polynomial-like functions,}. 

Consider a smooth function $f$ with  $n$ distinct
real zeros   
$x_{1}^{(0)} < x_{2}^{(0)} <\ldots <x_{n}^{(0)}$ in some interval $I\subseteq \bR$. Then, by Rolle's theorem, $f'$ has at
least $(n-1)$ zeros, $f''$ has at least $(n-2)$ zeros, ... , $f^{(n-1)}$
has at least one zero in the open interval $(x_{1}^{(0)},x_{n}^{(0)})$.
We are interested in  smooth functions $f$ with  $n$  real simple zeros in  $I$  such that  for all $i=1,\ldots,n,$ the
$i$th derivative  $f^{(i)}$ has exactly $n-i$ real simple 
zeros  in $I$  denoted by $x^{(i)}_1<x_2^{(i)}<...<x_{n-i}^{(i)}.$ Note, in particular, that $f^{(n)}$ is non-vanishing in $I$. 

\medskip 
\noindent
  \begin{definition} A smooth  function $f$ defined 
     in an interval $I$   is called {\bf polynomial-like of  degree $n$\/} 
       if $f^{(n)}$ does not  vanish in $I$. A polynomial-like function of degree $n$  in $I$ with $n$ simple real zeros   is called  {\bf real-rooted}. 
\end{definition}  

  By Rolle's theorem, $n$ is  the maximal possible number of real zeros in $I$ for a polynomial-like function of degree $n$.   An obvious example of a  real-rooted polynomial-like function of degree $n$  on  $\bR$ is a  usual  real polynomial of  degree $n$ 
with all real and distinct zeros.  Observe also that if a polynomial-like function $f$ of degree $n$ is real-rooted   in $I$, then  for all $i<n,$ 
its derivatives $f^{(i)}$ are also real-rooted of degree $n-i$ in the same interval. In the above notation the
following system of inequalities holds: 
\begin{equation}\label{eq:1}
         x^{(i)}_l<x^{(j)}_l<x^{(i)}_{l+j-i},\quad\text{    for  } i<j\le n-l. 
         \end{equation}
    We call \eqref{eq:1}  the  system of 
{\it  standard Rolle's restrictions.} (It is worth mentioning that the standard Rolle's restrictions define the well-known Gelfand-Tsetlin polytope, see e.g., \cite{deLoMcAl}.) 
 
With any real-rooted polynomial-like function $f$ of degree $n$, one can associate its configuration 
$\mathcal A_{f}$ of  $\binom {n+1} 2$ zeros $\{x^{(i)}_{l}\}$ of
$f^{(i)},$  for $i=0,\ldots,n-1$ and $1\le l\le n-i $,  by taking first all 
$x^{(0)}_{l}$, then all $x^{(1)}_{l},$ etc. 
\medskip

\medskip \noindent
\begin{problem}[V.~Kostov, B.~Sh., looks ugly, see \cite{ShShadd}]\label{pr:ug} What additional restrictions besides \eqref{eq:1} exist 
on configurations $\mathcal A_{f}=\{x^{(i)}_{l}\}$ coming from real-rooted polynomial-like  
 functions of a given degree $n$?  Or,  more ambitiously, given a configuration  $\mathcal A=\{x^{(i)}_{l}\vert\; 
i=0,\ldots,n-1;\,l=1,\ldots n-i\}$ of 
$\binom {n+1} 2$ real numbers satisfying  standard Rolle's 
restrictions,  is it possible to determine if there exists a real-rooted polynomial-like $f$ of degree $n$ such that 
$\mathcal A_{f}=\mathcal A$? 
\end{problem}

In the simplest non-trivial case $n=3,$ Problem~\ref{pr:ug}  was solved in \cite{ShShadd}, but the general case is widely open. 

\medskip
Substituting each zero of $p$ by the symbol $0$, each zero 
     of $p'$ by $1$, ... , the unique zero of $p^{(n-1)}$ by $(n-1)$ respectively, we get a {\it symbolic sequence of $p$} of 
length $\binom {n+1} {2}$ 
     with $n$ occurrences of $0$, $(n-1)$ occurrences of $1$,..., one occurrence of $(n-1)$. Standard Rolle's restrictions result in the condition 
 that between any two consecutive occurrences 
     of the symbol $i$ in  such a sequence, one has exactly one occurrence of the symbol $i+1$. 
 
     For example, there are only two possible symbolic sequences  $012010$ and $010210$,  for 
     $n=3$. For $n=4,$ there are 12 such 
     sequences $0123012010$, $0120312010$, $0120132010$, $0102312010$, 
     $0102132010$, 
     $0123010210$, $0120310210$, $0120130210$, $0120103210$, $0123010210$,       
     $0102130210$, $0102103210$. A patient reader will find that for $n=5,$ there are 
     $286$ such sequences. 
 
  If we denote by $\flat_{n}$ the number of all possible symbolic sequences 
  of length $n$,  then it  is possible to calculate this number  explicitly.  
It turns out to be equal to 
$$\flat_n=\binom{n+1} 
{2}!\frac{1!2!...(n-1)!}{1!3!...(2n-1)!}.$$  
%Since this calculation requires a separate treatment we refer an  interested reader to \cite  {Th}. 
 
%We can now formulate  a natural discrete analog of the main problem from the previous section which makes sense for  %usual polynomials.  

\medskip \noindent 
\begin{problem}[looks ugly] What symbolic sequences can occur for strictly real-rooted polynomials of degree $n$? 
%We  call such sequences  {\it realizable}. 
\end{problem}

\bigskip
\noindent
{\bf IX. Around "Hawaiian`` conjecture.} %See \cite{Tya}. 

\medskip
A simple observation that if a real polynomial $p(x)$  has all real and simple zeros  then the function $\frac
{p'(x)}{p(x)}$ is (locally) strictly monotone was apparently known to Gauss.
%In what follows we refer to  real polynomials with only (simple) real zeros
%as {\em (strictly) hyperbolic} and real polynomials with no real zeros as {\em
%(strictly) elliptic}. (An elliptic polynomial is necessarily of even degree.)
%Using this terminology 
We can reformulate the above observation in
the form of the classical Laguerre inequality:

\begin{lemma}
If $p(x)$ has only simple real zeros  then the polynomial
$P_{1}(x)=(p'(x))^2-p(x)p''(x)$ is strictly positive.
\end{lemma}

\medskip
Refinement of this observation constitutes  the  "Hawaiian`` conjecture saying that   for any real polynomial $p(x)$ with simple real zeros (and arbitrary complex-conjugate pairs of zeros), 
\begin{equation}
\sharp_{r}\left[(p'(x))^2-p(x)p''(x)\right] \le \sharp_{nr}p(x).
\label{eq:Haw}
\end{equation}
Here $\sharp_{r}q(x)$ (resp. $\sharp_{nr}q(x)$)
stands for the number of
real (resp.  non-real) zeros of a polynomial $q(x)$ with
real coefficients.

\medskip
"Hawaiian`` conjecture  was settled in \cite{Tya} by elementary but tedious calculations. 
%Here we try to find a proper context for conjecture (\ref{eq:Haw}) which
%leads to its natural generalizations as well as to some further developments.
Observe that   polynomial $P_{1}(x)=(p'(x))^2-p(x)p''(x)$
   appears not only as the numerator of $\left(\frac
   {p'(x)}{p(x)}\right)'$ but also in a different situation discovered
 by J.~Jensen, \cite {Je}  around 1910.  Namely, consider the
   function
   \begin{equation}
   \Phi_p(x,y)=\vert{p(x+iy)}\vert^2.
   \label{phi}
   \end{equation}
   $\Phi_p(x,y)$ is a real-analytic nonnegative function in $(x,y)$ whose
   zeros correspond to the zeros of $p(z)=p(x+iy)$. Assuming that
   $\deg p(x)=k$ and expanding
   $\Phi_p(x,y)$ in the variable $y,$ one gets
   \begin{equation}
   \Phi_p(x,y)=\sum_{i=0}^{k} P_{i}(x)\frac {y^{2i}}{(2i)!},\label{func}
   \end{equation}
   where
   \begin{equation}
   P_{i}(x)=\sum_{j=0}^i(-1)^{i+j}\binom {2i}{j}
   p^{(j)}(x)p^{(2i-j)}(x)\label{fi}.
   \end{equation}
   {In particular, $P_{0}(x)=p^2(x)$ and
   $P_{1}(x)=(p'(x))^2-p(x)p''(x)$.) The following  explicit formula
   \begin{equation}
   P_{i}(x)=p^2(x)\sum_{(l_{1},\ldots,l_{2i})}
    \frac{(2i)!}{(x-x_{l_{1}})^2\ldots (x-x_{l_{2i}})^2}\label{fie}
    \end{equation}
 is valid. Here $p(x)=(x-x_{1})(x-x_{2})\ldots(x-x_{k})$
    and the summation is taken over all $2i$-tuples (with repetitions
    in case of multiple zeros). %Note that $\deg P_{i}(x)=2(k-i)$. 
The latter formula immediately implies
    the following  criterion of real-rootedness known to J.~Jensen.

    \begin{proposition}[$1$-st criterion of  real-rootedness]\label{cr1}
A polynomial $p(x)$ of degree $k$  has simple real zeros if and only if  the polynomials $P_{i}(x)$ are
strictly positive for all $i=1,\ldots, k-1$.
\end{proposition}

Somewhat later G.~Polya studying a  number of unpublished notes left after Jensen's untimely death in 1912 discovered a different criterion of real-rootedness,% see \cite {PolJe}, p. 301..%

    \begin{proposition}[$2$-nd criterion of real-rootedness]\label{cr2}
A polynomial $p(x)$ of degree $k$ has real simpe zeros if and only if  the polynomials
\begin{equation}\label{weighted}
G_{i}(x)=(k-i)(p^{i}(x))^2-(k-i+1)p^{(i-1)}(x)p^{(i+1)}(x)
\end{equation}  are strictly  positive for all $i=1,\ldots, k-1$.
\end{proposition}

Propositions~\ref{cr1}-\ref{cr2} combined with the original "Hawaiian`` conjecture   motivate our questions of the number of real zeros for the families  $\{P_{i}(x)\}$
and $\{G_i(x)\}$ presented below. They are  based on extensive experiments with real polynomials of degree up to $6$.

\begin{conjecture}[B.~Sh., seems good] \label{conjweight}
For any real polynomial $p(x)$ of degree $k$ with simple real zeros, 
\begin{equation}
\sharp_{r}\left[(k-1)(p'(x))^2-kp(x)p''(x)\right] \le \sharp_{nr}p(x),
\end{equation}
i.e. "Hawaiian`` conjecture holds for $G_1(x)$ as well.
\end{conjecture}

\begin{corollary}[Conjectural] For any real polynomial $p(x)$ of degree $k$  with simple real
zeros,  
\begin{equation}
\sharp_{r}G_{i}(x) \le \min\{{\deg{G_i(x)},\sharp_{nr}p(x)}\}.
\label{eq:leiHaw}
\end{equation}
%where (as above) $\sharp_{r}q(x)$ resp. $\sharp_{nr}q(x)$
%stands for the number of
%real resp.  nonreal zeros of a considered polynomial $q(x)$ with
%real coefficients.
    \label{conj1}
    \end{corollary}

    Additionally, we claim
    \begin{conjecture}[B.~Sh]\label{conjwplus} For any real polynomial $p(x)$ of even degree, 
    \begin{equation}
\sharp_{r}\left[(k-1)(p'(x))^2-kp(x)p''(x)\right] + \sharp_{r}p(x)>0,
\end{equation}
    \end{conjecture}

    The latter inequality is trivially satisfied for odd degree polynomials.

\begin{conjecture}[B.~Sh] For any degree $k$ polynomial $p(x)$ with real
coefficients,  
\begin{equation}
\sharp_{r}P_{i}(x) \le \min\{{\deg{P_i(x)},k}\}.
\label{eq:HighLag}
\end{equation}
    \label{conj2}
    \end{conjecture}


\begin{thebibliography}{8}


\bibitem[BKSh]{BKSh} P.~Br\"and\'en, I.~Krasikov, B.~Shapiro, Elements of P\'olya-Schur theory in finite diffrence setting, arXiv:1204.2963.

\bibitem[CLR]{CLR} M-D.~Choi, T-Y.~Lam, B.~Reznick, Real zeros of positive definite forms. I, Math. Z., vol. 171 (1980), 1--26. 

\bibitem[deLoMcAl]{deLoMcAl} J.~A.~De Loera, T.~B.~McAllister, {Vertices of Gelfand-Tsetlin Polytopes}, arXiv:math/0309329.

\bibitem [GNSh]{GNSh} A.~Gabrielov, D.~Novikov, and B.~Shapiro, Mystery of point charges, Proc. London Math. Soc. (3) vol. 95(2) (2007), 443--472. 

\bibitem[Fi]{Fi}  S.~Fisk,  Polynomials, roots, and interlacing, arXiv:math/0612833. 

\bibitem[FKSh]{FKS} J.~Forsg\aa rd, Vl.~Kostov, B.~Shapiro, {Could Ren\'e Descartes have known this?}, arXiv:1501.00856.   

\bibitem[FSh]{FSh} J.~Forsg\aa rd, B.~Shapiro, {Discriminant amoebas,  complex zero decreasing sequences, and the Karlin problem}, in preparation.

\bibitem[Gr]{Gr}  D.~J.~Grabiner, Descartes' Rule of Signs: Another Construction, The American Mathematical Monthly vol. 106 (1999), 854--856.

\bibitem[Je]{Je} J.~L.~W.~V.~Jensen, Recherches sur la th\'eorie des \'equations, Acta Math., vol. 36 (1913), 181--195.   

\bibitem[Kh]{Kh} A.~G.~Khovanskii, Fewnomials,  AMS translated monographs, USA, V.88, viii+139, 1991.

\bibitem[Ki]{Ki} K.~Killian, A remark on Maxwell's conjecture for planar charges, Complex Variables and Elliptic Equations, vol. 54(12) (2009), 1073--1078.

\bibitem[KST]{KST} V.~Kostov, B.~Shapiro, and M.~Tyaglov, Maximal univalent disks of real rational functions and Hermite-Biehler polynomials,  Proc. Amer. Math. Soc. vol. 139(5), (2011) 1625--1635.

\bibitem[La]{La}R.~E.~Langer, On the zeros of exponential sums and integrals, Bull. Amer. Math. Soc. vol. 37 (1931), 213--239.


%\bibitem [MS] {MS}  G.~M\'asson and B.~Shapiro,  A note on  polynomial eigenfunctions of a hypergeometric type
%operator,  Experiment. Math. vol. 10, issue 4 (2001) 609--618. 

\bibitem[Ma]{Maxwell} J.~C.~Maxwell,
      {\sl A Treatise on Electricity and Magnetism}, {\bf vol. 1},
      Republication of the 3rd revised
      edition, Dover Publ. Inc., 1954.

\bibitem[OlPe]{OlPe} O.~A.~Oleinik and I.~G.~Petrovskii, On the topology of real algebraic hypersurfaces (Russian),
Izv. Acad. Nauk SSSR vol. 13 (1949) 389--402. English transl.: Amer. Math. Soc. Transl. vol. 7
(1962) 399--417.

\bibitem[Pe]{Pe} R.~Peretz, Application of the argument principle to Maxwell's Conjecture for three point charges, 
 Complex Variables and Elliptic Equations , vol. 58(5) (2013),  1--11. 

\bibitem [ShSh] {ShShadd}   B.~Shapiro,  M.~Shapiro,  A few riddles behind Rolle's theorem,   Amer. Math. Monthly, vol.  119 (9)  (2012), 787--793. 

\bibitem[Sch] {Sch} W. M. Schmidt, The zero multiplicity of linear recurrence sequences, Acta Math. vol. 182 (1999), 243--282.

\bibitem [SSm]{SSm} T.~Sheil-Small,  {Complex polynomials}, Cambridge Studies
in Adv. Math. vol. 75, Cambridge Univ. Press, Cambridge, UK, (2002).

\bibitem[Tya]{Tya} M.~Tyaglov,  On the number of real critical points of logarithmic derivatives and the Hawaii conjecture. J. Anal. Math. vol. 114 (2011), 1--62.

\bibitem[Wa]{Wa} Yi.~Wang, Equilibrium Points of Potential Fields Produced by Positive Point Charges, ccs.math.ucsb.edu/senior-thesis/YifeiWang.pdf

\end{thebibliography}
\end{document}